\documentclass[11pt]{amsart}

\usepackage{mathrsfs}
\usepackage{extarrows}
\usepackage{amssymb, amsmath, amsthm}
\usepackage{graphicx,verbatim}
\usepackage{rotating}
\usepackage[OT2,T1]{fontenc}
\usepackage[all,cmtip]{xy}

    \let\oldmarginpar\marginpar
    \renewcommand\marginpar[1]{\-\oldmarginpar[\raggedleft\footnotesize #1]{\raggedright\footnotesize #1}}

    \theoremstyle{plain}
    
    \newtheorem{theorem}{Theorem}[section]
    \newtheorem{proposition}[theorem]{Proposition}
    \newtheorem{corollary}[theorem]{Corollary}
    \newtheorem{lemma}[theorem]{Lemma}
    \newtheorem{conjecture}[theorem]{Conjecture}
    \newtheorem*{theorem*}{Theorem}
    \newtheorem*{proposition*}{Proposition}
    \newtheorem*{corollary*}{Corollary}
    \newtheorem*{lemma*}{Lemma}
    \newtheorem*{conjecture*}{Conjecture}

    \theoremstyle{definition}
    \newtheorem{definition}[theorem]{Definition}
    \newtheorem{example}[theorem]{Example}
    \newtheorem*{definition*}{Definition}
    \newtheorem*{example*}{Example}

    \theoremstyle{remark}
    \newtheorem{remark}[theorem]{Remark}
    \newtheorem*{remark*}{Remark}

    \newenvironment{myproof}[1][\noindent\textit{\sc \textbf{Proof.} }]
    {#1}{$\hfill\qedsymbol\newline$}

    \providecommand{\n}{$\newline$}

    \newcommand{\s}{\sigma}

    \newcommand{\bZ}{\mathbb{Z}}
    \newcommand{\bQ}{\mathbb{Q}}
    
    \newcommand{\bC}{\mathbb{C}}

    \newcommand{\bP}{\mathbb{P}}

    \newcommand{\cC}{\mathcal{C}}

    \newcommand{\cM}{\mathcal{M}}

    \newcommand{\cX}{\mathcal{X}}

    \newcommand{\sE}{\mathscr{E}}
    \newcommand{\sF}{\mathscr{F}}

    \newcommand{\sI}{\mathscr{I}}

    \newcommand{\sL}{\mathscr{L}}
    \newcommand{\sM}{\mathscr{M}}
    \newcommand{\sN}{\mathscr{N}}
    \newcommand{\sO}{\mathscr{O}}
    
    \newcommand{\sQ}{\mathscr{Q}}

    \newcommand{\sT}{\mathscr{T}}

    \newcommand{\rk}{\operatorname{rk }}

    \newcommand{\Hom}{\operatorname{Hom}}
    
    \newcommand{\et}{\operatorname{et}}
    
    \newcommand{\pr}{\operatorname{pr}}
    \newcommand{\codim}{\operatorname{codim}}
    \newcommand{\Alb}{\operatorname{Alb}}
    
    \newcommand{\Spec}{\operatorname{Spec}}
    \newcommand{\ev}{\operatorname{ev}}

    \newcommand{\Bl}{\operatorname{Bl}}

    \newcommand{\arrow}{\rightarrow}

    \renewcommand\bar\overline
    \renewcommand\tilde\widetilde

    \newcommand{\eqn}[1]{\begin{eqnarray*}#1\end{eqnarray*}}
    \newcommand{\eqnn}[1]{\begin{eqnarray}#1\end{eqnarray}} 

    \newcommand{\theo}[1]{\begin{theorem}#1\end{theorem}}
    \newcommand{\theon}[1]{\begin{theorem*}#1\end{theorem*}}

    \newcommand{\prop}[1]{\begin{proposition}#1\end{proposition}}

    \newcommand{\lem}[1]{\begin{lemma}#1\end{lemma}}

    \newcommand{\cor}[1]{\begin{corollary}#1\end{corollary}}

    \providecommand{\defn}[1]{\begin{definition}#1\end{definition}}

    \newcommand{\rem}[1]{\begin{remark}#1\end{remark}}

    \renewcommand{\proof}[1]{\begin{myproof} #1\end{myproof}}

\renewcommand{\qedsymbol}{\ensuremath{\Box}}

\usepackage{hyperref}
\hypersetup{
    colorlinks=true,       
    linkcolor=blue,          
    citecolor=blue,        
    filecolor=blue,      
    urlcolor=blue           
}

\numberwithin{equation}{section}

\theoremstyle{remark}

\newtheorem*{acknowledgements}{acknowledgements}

\address{Frank Gounelas, Institut f\"ur Mathematik, Humboldt Universit\"at zu Berlin,
Unter den Linden 6, 10099 Berlin.}
\author{Frank Gounelas}
\title{Free curves on varieties}
\date{\today}

\begin{document}
\maketitle
\begin{abstract}
    We study various generalisations of rationally connected varieties, allowing the
    connecting curves to be of higher genus. The main focus will be on free curves $f:C\to
    X$ with large unobstructed deformation space as originally defined by Koll\'ar, but we
    also give definitions and basic properties of varieties $X$ covered by a family of
    curves of a fixed genus $g$ so that through any two general points of $X$ there passes
    the image of a curve in the family. We prove that the existence of a free curve of
    genus $g\geq1$ implies the variety is rationally connected in characteristic zero and
    initiate a study of the problem in positive characteristic.
\end{abstract}

\section{Introduction}

    Let $k$ be an algebraically closed field. A smooth projective rationally connected
    variety, originally defined in \cite{campana92} and \cite{kmm92}, is a variety such
    that through every two general points there passes the image of a rational curve. In
    characteristic zero this is equivalent to the notion of a separably rationally
    connected variety, given by the existence of a rational curve $f:\bP^1\to X$ such that
    $f^*\sT_X$ is ample. In characteristic $p$, however, one has to distinguish between these
    two notions. Deformations of a morphism $f:\bP^1\to X$ are controlled by the sheaf
    $f^*\sT_X$, hence studying positivity conditions of this bundle is intimately tied to
    deformation theory and the existence of many rational curves on $X$. Rationally
    connected varieties have especially nice properties and an introduction to the theory
    is contained in \cite{kollar} and \cite{debarre}. Note in particular the important
    theorem of Graber-Harris-Starr \cite{ghs} (and de Jong-Starr \cite{djs} in positive
    characteristic) which we will make repeated use of throughout this paper, which says
    that a separably rationally connected fibration over a curve admits a section. An
    equivalent statement in characteristic zero is that the maximal rationally connected
    (MRC) quotient $R(X)$ is not uniruled (see \cite[IV$.5.6.3$]{kollar}), although this can
    fail in positive characteristic. 
    
    \n
    In this paper we study various ways in which a variety can be connected by higher
    genus curves. After an introductory section with auxiliary results on vector
    bundles on curves and Frobenius, we consider first varieties which admit a morphism
    from a family of curves of fixed arithmetic genus $g$ whose product with itself
    dominates the product of the variety with itself and call these varieties ``genus $g$
    connected'', generalising the notion of there being a rational curve through two
    general points. We also consider $C$-connected varieties, where there exists a family
    $C\times U \to X$ of a single smooth genus $g$ curve $C$ such that $C\times C\times
    U\to X\times X$ is dominant. Mori's Bend and Break result allows us to produce
    rational curves going through a fixed point given a higher genus curve which has large
    enough deformation space. For example, in Proposition \ref{cconisuniruled} as an easy corollary, we show that over
    any characteristic, if for any two general points of a smooth projective variety $X$ with $\dim X\geq3$ there passes
    the image of a morphism from a fixed curve $C$ of genus $g$, then $X$ is uniruled. This fails for surfaces, where an
    example is provided.
    
    \n
    A stronger condition than the aforementioned is the existence of a morphism from a
    curve which deforms a lot without obstructions, as discussed for separably rationally
    connected varieties above. Namely, for $f:C\to X$ a morphism to a variety $X$ where
    $C$ is of any genus $g$, Koll\'ar \cite{kollar} defines $f$ to be free if $f^*\sT_X$
    is globally generated as a vector bundle on $C$ and also $H^1(C, f^*\sT_X)=0$. In the
    case of genus $g=0$ one must distinguish between free and very free curves.
    Geometrically, the former implies that $f:\bP^1\to X$ deforms so that its image covers all points
    in $X$ (hence $X$ is uniruled) whereas the latter that it can do so even fixing a point $x\in X$
    ($X$ rationally connected). If $g\geq1$, however, after defining an
    $r$-free curve to be one which deforms keeping any $r$ points fixed, we show that the
    notions of the existence of a free ($0$-free) and very free ($1$-free) curve coincide
    and in fact are equivalent with the existence of a curve $f:C\to X$ such that
    $f^*\sT_X$ is ample.

    \theon{(see \ref{rfreeiffample})
        Let $X$ be a smooth projective variety and $C$ a smooth projective curve of genus
        $g\geq1$ over an algebraically closed field $k$. Then for any $r\geq0$, there
        exists an $f:C\to X$ which is $r$-free if and only if there exists a morphism
        $f':C\to X$ such that $f'^*\sT_X$ is ample.
    }

    Work of Bogomolov-McQuillan (see \cite{bm}, \cite{ksct}) on foliations which restrict
    to an ample bundle on a smooth curve sitting inside a complex variety $X$ shows that
    the leaves of such a foliation are not only algebraic but in fact have rationally
    connected closures. From the above, one deduces this result in the case of the
    foliation $\sF=\sT_X$, complementing the currently known connections between existence
    of curves with large deformation space and rationally connected varieties (cf. the
    uniruledness criterion of Miyaoka \cite{miyaoka}). Our proof emphasises the use of
    free curves and $C$-connected varieties, in particular with a view towards similar
    results in positive characteristic. 
    
    \theon{(see \ref{mrctowerispoint})
        \label{maintheointro}
        Let $X$ be a smooth projective variety over an algebraically closed field of
        characteristic zero and let $f:C\to X$ be a smooth projective curve of genus
        $g\geq1$ such that $f^*\sT_X$ is globally generated and $H^1(C, f^*\sT_X)=0$. Then
        $X$ is rationally connected.
    }

    In the sixth section we study the particular case of elliptically connected varieties
    (i.e.\ genus one connected varieties) where, even allowing a covering family of genus
    $1$ curves to vary in moduli, one can prove the following theorem.

    \theon{(Theorem \ref{elconclassification})
        Let $X$ be a smooth projective variety over an algebraically closed field of
        characteristic zero. Then the following two statements are equivalent
        \begin{enumerate}
            \item There exists $\cC\to U$ a flat projective family of irreducible genus
                $1$ curves with a map $\cC\to X$ such that $\cC\times_U\cC\to X\times X$ is
                dominant.
            \item $X$ is either rationally connected or a rationally connected fibration
                over a curve of genus one.
        \end{enumerate}
    }

    In positive characteristic, at this point we have not been able to prove that the
    existence of a higher genus free curve implies the existence of a very free rational
    curve (which would mean that $X$ is separably rationally connected). We work however
    in this direction, establishing this result in dimensions two (with a short discussion
    about dimension three) and furthermore by studying algebraic implications of the
    existence of a free higher genus curve, such as the vanishing of pluricanonical forms
    and triviality of the Albanese variety. In the final section we give an example of a
    threefold in characteristic $p$ whose MRC quotient is rationally connected and which
    has infinite fundamental group.

    \n
    The study of rational curves on varieties is an important and active area of research, and shedding light on the
    existence of rational curves coming from the deformation theory of higher genus curves is a theme explored in a
    variety of sources, for example the minimal model program or \cite{bdppnew}. Aside from the unresolved difficulties
    arising in positive characteristic, the author expects uniruledness and rational connected results of the type
    presented in this article to be of use in moduli theory.

    \acknowledgements{
        The contents of this paper are from the author's thesis under the supervision of Victor Flynn, whom I would like
        to thank for his continuous encouragement. I am indebted to Damiano Testa for the many hours spent helping with
        the material of this paper and to Johan de Jong not only for the hospitality at Columbia University but also for
        helping improve the contents of this paper. I would also like to thank Jason Starr and Yongqi Liang for
        comments, J\'anos Koll\'ar for pointing out a similar construction to that in the last section and Mike Roth for
        showing me how abelian surfaces are $C$-connected. The anonymous referee's numerous suggestions and corrections
        also significantly improved this paper. This research was completed under the support of EPSRC grant number
        EP/F$060661/1$ at the University of Oxford. 
    }

\section{Ample vector bundles and Frobenius}

    We begin with some results concerning positivity of vector bundles on curves. Recall
    that a locally free sheaf $\sE$ on a scheme $X$ is called ample if $\sO_{\bP(\sE)}(1)$
    has this property. Equivalent definitions involving global generation of $\sF\otimes
    S^n(\sE)$ for $\sF$ a coherent sheaf and $n$ large enough, and also cohomological
    vanishing criteria can be found in \cite{hartample}. Ampleness on curves can be
    checked using various criteria such as the following.

    \lem{
        Let $C$ be a smooth projective curve of genus $g\geq2$ over an algebraically
        closed field of characteristic zero and $\sE$ a locally free sheaf on $C$ such that
        $H^1(C, \sE)=0$. It follows that $\sE$ is ample.
        \label{fujitalemma}
    }
    \proof{
        From \cite[Theorem $2.4$]{hartcurveample}, it suffices to show that every
        non-trivial quotient locally free sheaf of $\sE$ has positive degree. Let $\sE\to
        \sE'\to0$ be a quotient. From the long exact sequence in cohomology we see that
        $H^1(C,\sE')$ is also $0$. From the Riemann-Roch formula $\deg \sE' = h^0(C, \sE') 
        + (\rk \sE')(g-1)$ and since $g\geq2$ we deduce that $\deg \sE' >0$.
    }

    Note that Hartshorne's ampleness criterion only works in characteristic zero. More
    generally, over any characteristic if we further assume that our locally free sheaf is
    globally generated then the same result holds so long as the genus is at least one.

    \prop{
        Let $C$ be a smooth projective curve of genus $g\geq1$ over an algebraically closed
        field $k$ and $\sE$ a globally generated locally free sheaf on $C$ such that $H^1(C,
        \sE)=0$. Then $\sE$ is ample.
        \label{ggandh1impliesample}
    }
    \proof{
        Since $\sE$ is globally generated, there exists a positive integer $n$ such that
        $\sO_C^{\oplus n} \to\sE\to 0$ is exact. This gives (see \cite[ex.
        II$.3.12$]{hartshorne}) a closed immersion of the respective projective bundles
        $\bP(\sE) \hookrightarrow \bP^{n-1}_C$. By projecting onto the first factor we
        have the following diagram
        \eqn{
            \xymatrix{
                \bP(\sE) \ar@{^{(}->}^-i[r]\ar[dr]^-{\pi} & \bP^{n-1}\times
                C\ar[d]^{\pr_2}\ar[r]^-{\pr_1} & \bP^{n-1} \\ & C &
            }
        }
        and from \cite[II$.5.12$]{hartshorne} we have $\pr_1^*\sO_{\bP^{n-1}}(1) =
        \sO_{\bP^{n-1}_C}(1)$. Also, since $i$ is a closed immersion it follows that
        $i^*\sO_{\bP^{n-1}_C}(1) = \sO_{\bP^{n-1}_C}(1)|_{\bP(\sE)} = \sO_{\bP(\sE)}(1)$
        which concludes that $i^*\pr_1^*\sO_{\bP^{n-1}}(1) = \sO_{\bP(\sE)}(1)$. To show
        that $\sE$ is an ample locally free sheaf on $C$ it is enough to show that this
        invertible sheaf is ample. Since we know that $\sO_{\bP^{n-1}} (1)$ is ample
        though, it is sufficient to show that $i\circ\pr_1$ is a finite morphism. Since it
        is projective, we need only show that it is quasi-finite. Hence assuming that the
        fibre of $i\circ\pr_1$ over a general point $p\in\bP^{n-1}$ is not finite, it must
        be the whole of $C$. We now embed this fibre $j:C\to\bP(\sE)$ as a section to
        $\pi$ and pull back the surjection $\pi^*\sE \to \sO_{\bP(\sE)}(1)$ via $j$,
        obtaining $j^*\sO_{\bP(\sE)}(1)$ as a quotient of $j^*\pi^*\sE = \sE$ (see
        \cite[II.$7.12$]{hartshorne}). However $\pr_1\circ i\circ j:C\to \bP^{n-1}$ is a
        constant map so $j^*\sO_{\bP(\sE)}(1) = \sO_C$. Taking cohomology of the
        corresponding short exact sequence given by this quotient, we obtain a
        contradiction since $H^1(C,\sE) = 0$ whereas $H^1(C,\sO_C)$ is not trivial for
        $g\geq1$.
    }

    In Proposition \ref{ampleisvgincharp} below we will prove that given an ample bundle
    on a curve in positive characteristic, then after pulling back by Frobenius, we can
    make this bundle be globally generated and have vanishing first cohomology.

    \lem{
        Let $C$ be a smooth projective curve over an algebraically closed field $k$,
        $d\geq0$ an integer and $\sE$ a locally free sheaf on $C$. If $H^1(C, \sE(-D))=0$
        for all effective divisors $D$ of fixed degree $d$ then for $d'<d$ it follows
        that $H^1(C, \sE(-D'))=0$ and $\sE(-D')$ is globally generated for all
        effective divisors $D'$ of degree $d'$.
        \label{appendixvfimpliesf}
    }
    \proof{
        The first result follows from the short exact sequence
        \eqn{
            0\to \sE(-D'-R) \to \sE(-D') \to \sE(-D')|_R \to 0
        }
        where $R$ is an effective divisor of degree $d-d'$. For the second, let $p\in C$.
        From the first part we have $H^1(C,\sE(-D'-p))=0$ since $D'+p$ is an effective
        divisor of degree $d'+1\leq d$ so the following sequence is exact
        \eqn{
            0\to H^0(C,\sE(-D'-p)) \to H^0(C,\sE(-D')) \to \sE(-D')\otimes k(p) \to 0.
        }
        Hence $\sE(-D')$ is globally generated at $p$ and the result follows.
    }

    A partial converse to Proposition \ref{ggandh1impliesample} in characteristic $p$ is
    given in \cite[Proposition $9$]{ksct}, using $\bQ$-twisted vector bundles as in
    \cite[II$.6.4$]{lazarsfeld2}. We prove the following different version of this result. 

    \prop{
        Let $C$ be a smooth projective curve of genus $g$ over an algebraically closed
        field $k$ of characteristic $p$ and let $\sE$ be an ample locally free sheaf on
        $C$. Let $B\subset C$ be a closed subscheme of length $b$ and ideal sheaf $\sI_B$.
        Then there exists a positive integer $n$ such that $H^1(C^{(n)}, F_n^*\sE\otimes
        \sI_B)=0$ and $F_n^*\sE\otimes \sI_B$ is globally generated on $C^{(n)}$ where
        $F_n:C^{(n)}\to C$ the $n$-fold composition of the $k$-linear Frobenius morphism.
        \label{ampleisvgincharp}
    }
    \proof{
        We proceed by induction. First, assume we can write $\sE$ as an extension
        \eqn{
            0\to \sM\to \sE\to \sQ\to 0
        }
        where $\sM$ is an ample line bundle. If $\sQ$ is not torsion free, consider the
        saturation of $\sM$ in $\sE$ instead and take $\sQ$ as that quotient. Since $\sE$
        is ample, so is its quotient $\sQ$. Note also that the rank of $\sQ$ is one less
        than that of $\sE$ and that if we can prove the result for $\sQ$ then we will have
        it for $\sE$ too by considering cohomology of the appropriate exact sequences. We
        thus reduce to the case of $\sE=\sL$ an invertible sheaf of positive degree (since
        it is ample). An invertible sheaf $\sL$ pulls back under the $n$-fold composition
        of the linear Frobenius morphism to an invertible sheaf $F_n^*\sL$ of degree
        $p^n\deg \sL$. To show that $H^1(C^{(n)}, F_n^*\sL\otimes \sI_B)=0$, it is
        equivalent by Serre duality to show that $\Hom_{C^{(n)}}(F_n^*\sL,
        \sO_{C^{(n)}}(B)\otimes \omega_{C^{(n)}})=0$. Since the invertible sheaf
        $\sO_{C^{(n)}}(B)\otimes \omega_{C^{(n)}}$ has degree $b+2g-2$ and by picking $n$
        large enough, we can ensure $p^n\deg \sL > b+2g-2$ from which we obtain
        $H^1(C^{(n)}, F_n^*\sL\otimes \sI_B)=0$ and hence $H^1(C^{(n)}, F_n^*\sE\otimes
        \sI_B)=0$ for a locally free sheaf of any rank.

        \n
        To show that $F_n^*\sE\otimes \sI_B$ is globally generated, pick a point $q\in C$.
        Then $\sI_B\otimes \sI_q$ has length $b+1$ and from the discussion above $H^1(C^{(n)},
        F_n^*\sE\otimes \sI_B\otimes \sI_q)$ vanishes when $p^n\deg L > b+1+2g-2$ so we can just
        pick $n$ large enough to fit this condition. Now, by taking the long exact
        sequence in cohomology of 
        \eqn{
            0\to F_n^*\sE\otimes \sI_B\otimes \sI_q \to F_n^*\sE\otimes \sI_B \to (F_n^*\sE\otimes
            \sI_B)\otimes k(q)\to 0
        }
        we conclude that $F_n^*\sE\otimes \sI_B$ is globally generated.

        \n
        That $\sE$ can not be written as an extension of $\sM$ an ample line bundle and a
        quotient locally free sheaf $\sQ$ is equivalent to $H^0(C, \sE\otimes\sM^{-1})=0$.
        However there exists a positive integer $m$ and an ample line bundle
        $\sM_{C^{(m)}}$ on $C^{(m)}$ for which $H^0(C^{(m)}, (F_m^*\sE) \otimes
        \sM_{C^{(m)}}^{-1}) \neq 0$ and we proceed as before with the sheaf $(F_m^*\sE)$.
    }
    
\section{Definition of curve connectedness: Covering families}

    We now define various ways in which a variety can be covered by curves, generalising
    the notion of a rationally connected varieties (see \cite[IV]{kollar}).

    \defn{
        We say that a variety $X$ over a field $k$ is \textit{connected by genus $g\geq0$
        curves} (resp. \textit{chain connected by genus $g$ curves}) if there exists a
        proper flat morphism $\cC\arrow Y$, for a variety $Y$, whose geometric fibres are
        irreducible genus $g$ curves (resp. connected genus $g$ curves) such that there is
        a morphism $u: \cC \rightarrow X$ making the induced morphism $u^{(2)}:
        \cC\times_Y \cC \rightarrow X \times_k X$ dominant. \label{defngcon}
    }

    We say $X$ is separably (chain) connected by genus $g$ curves if $u^{(2)}$ is smooth
    at the generic point. Note that the notion of separability is redundant in
    characteristic zero due to generic smoothness. A genus zero connected variety is
    rationally connected. A variety which is connected by genus one curves will be called
    (with a slight abuse of notation) elliptically connected. The relevant moduli spaces
    which we will be considering are the following. Let $\pi:\cC\to S$ be a flat
    projective curve over an irreducible scheme $S$ and let $B\subset \cC$ be a closed
    subscheme that is flat and finite over $S$. Let $p:X\to S$ be a smooth
    quasi-projective scheme and $g:B\to X$ an $S$-morphism. The space (see
    \cite[II$.1.5$]{kollar} and \cite{mori}) $\Hom_S(\cC, X, g)$ parametrises
    $S$-morphisms from $\cC$ to $X$ keeping the points given by $g$ fixed. Restricting to
    the case where $S$ is the spectrum of an algebraically closed field $k$ we fix some
    notation of the following evaluation morphisms to be used in later sections
    \eqn{
        F:C\times\Hom(C,X,g)&\arrow& X \\
        \phi(p,f): H^0(C, f^*\sT_X\otimes \sI_B) &\arrow& f^*\sT_X\otimes k(p)
    }
    and similarly the double evaluation morphisms $F^{(2)}$ and $\phi^{(2)}(p,q,f)$ as in \cite[II$.3.3$]{kollar}.
    Secondly we consider the relative moduli space of genus $g$ degree $d$ stable curves with base point $t:P\to X$,
    denoted by $\bar{\cM}_g(X/S, d, t)$ as in \cite{ak} (originally \cite{fp}). By Bertini, we can always find a genus
    $g$ such that a projective $X$ is genus $g$ connected, the minimal such $g$ however is an interesting invariant of
    the variety. Finding higher genus covering families is an easy operation.

    \lem{
        Let $X$ be a genus $g$ (chain) connected smooth projective variety over an
        algebraically closed field $k$. Then if $g'\geq 2g-1$, $X$ is also genus $g'$
        (chain) connected.
        \label{genusgconnisgconforbiggerg}
    }
    \proof{
        Let $\cC/Y\to X$ be a family making $X$ a genus $g$ (chain) connected variety.
        From \cite[Theorem $50$]{ak} we have a projective algebraic space
        $Y'=\bar{\cM}_g'(\cC/Y,d)$ of finite type over $Y$ parametrising stable families
        of degree $d$ curves of genus $g'$ over $\cC\to Y$. The condition $g'\geq2g-1$
        coming from the Riemann-Hurwitz formula ensures that this moduli space is
        non-empty. From \cite[$12.9.2$]{acgh} there exists a normal scheme $Z$ finite and
        surjective over $Y'$ and a flat and proper family $\cX\to Z$ of stable genus $g$
        curves of degree $d$. Restricting to a suitable open subset $W\subset Z$
        parametrising irreducible curves we compose the family $\cX|_W\to W$ with the
        evaluation morphism to $X$ and the result follows.
    }

    An example of an elliptically connected variety over a non-algebraically closed field
    is given after the proof of Theorem \ref{elconclassification}. A much stronger
    condition is the existence of a family of curves which is constant in moduli.

    \defn{
        We say that a variety $X$ over a field $k$ is \textit{$C$-connected} for a curve
        $C$ if there exists a variety $Y$ and a map $u:C\times Y\to X$ such that the
        induced map $u^{(2)}: C\times C\times Y\to X\times X$ is dominant. If $u^{(2)}$ is
        also smooth at the generic point, then we say that $X$ is \textit{separably
        $C$-connected}.
        \label{defncconn}
    }

    Projective space is $C$-connected for every smooth projective curve $C$ whereas an
    example of a $C$-connected variety which is not rationally connected is $C \times
    \bP^n$ where $g(C)\geq1$. To see this let $(c_1,x_1),(c_2,x_2)$ be any two points in
    $C \times \bP^n$ and let $f:C\to \bP^n$ a morphism which sends $c_i\mapsto x_i$.
    Considering the graph of $f$ in $C\times\bP^n$ we have found a curve isomorphic to $C$
    which goes through our two points. Using parts $(3)$ and $(4)$ from Lemma
    \ref{gconproperties} below, the result follows. More generally, examples can also be
    constructed from Proposition \ref{rcfibrationovercurveiscconnected} below. The
    following are mostly straight forward generalisations of various results in
    \cite[IV$.3$]{kollar}.
    \lem{
        \label{gconproperties}
        The following statements hold for a variety $X$ over a field $k$ and $C$ a smooth
        projective curve.
        \begin{enumerate}
            \item If $X$ is genus $g$ connected and $X\dashrightarrow Y$ a dominant
                rational map to a proper variety $Y$, then $Y$ is also genus $g$
                connected. The same holds if $X$ is $C$-connected.
            \item A variety $X$ is $C$-connected if and only if there is a variety $W$,
                closed in $\Hom(C,X)$ such that $u^{(2)}: C\times C\times W\to X\times X$
                is dominant.
            \item If $X$ is defined over a field $k$ and $K/k$ is an extension of fields,
                then $X_K := X\times_k K$ is $C$-connected if and only if $X_k$ is.
            \item A variety $X$ over an uncountable algebraically closed field is
                $C$-connected if and only if for all very general $x_1,x_2\in X$ there
                exists a morphism $C\to X$ which passes through $x_1,x_2$.
            \item A variety $X$ over an uncountable algebraically closed field is genus
                $g$ connected if and only if for all very general $x_1,x_2 \in X$ there
                exists a smooth irreducible genus $g$ curve containing them. 
            \item Being rationally or elliptically connected is closed under connected
                finite \'etale covers of varieties.
        \end{enumerate}
    }
    \proof{
        To prove $(1)$, let $u:\cC/M\to X$ be a family making $X$ genus $g$ connected and
        denote by $u':\cC/M\dashrightarrow Y$ the composition. Restricting $u'$ to the
        generic fibre $\cC_{k(M)}$ we have a rational map $\phi: \cC_{k(M)}
        \dashrightarrow Y$. Since $Y$ is proper, by the valuative criterion
        of properness we can extend $\phi$ to a morphism
        $\phi: \cC_{k(M)} \to Y$. By spreading out to an open subset $M'\subseteq M$ (see
        \cite[IV$_3$ $8.10.5$]{ega} for properness and $11.2.6$ for flatness of the
        family) we obtain a family $\cC|_{M'}\to M'$ which makes $Y$ also genus $g$
        connected.
        
        \n
        Since being $C$-connected or connected by genus $g$ curves is a birational
        property, we may assume by compactifying that $X$ is projective. For $(2)$,
        consider $\Hom(C,X) = \cup R_i$ the decomposition into irreducible components. One
        direction of the statement is obvious, whereas for the other let $C\times W\to W$
        be a family which makes $X$ a $C$-connected variety. If $u_i:C\times R_i\to X$ is
        the evaluation morphism, then for some $i$ there is a morphism $h:W \to R_i$ such
        that $h(w) = [C_w\to X]$ for general $w\in W$. This implies that $u_i^{(2)}:
        C\times C\times R_i \to X\times X$ is also dominant. For one direction of $(3)$,
        pullback by $\Spec K\to \Spec k$. For the other, if $X_K$ is $C_K$-connected then
        from $(2)$ there is a positive integer $d$ such that the evaluation morphism
        $\ev_K^d: C_K\times C_K\times\Hom_d(C_K, X_K) \to X_K \times X_K$ is dominant.
        Because of the universal property of the $\Hom$-scheme, we have that $\Hom(C,X)
        \times_k K = \Hom(C_K, X_K)$ and $(\ev^d)_K = \ev_K^d$ so $ev^d$ is also dominant.
        
        \n
        If through every two very general points there passes the image of $C$ under some
        morphism, then the map $u^{(2)}: C\times C\times\Hom(C,X)\to X\times X$ is
        dominant. Since $\Hom(C,X)$ has at most countably many irreducible components the
        restriction of $u^{(2)}$ to at least one of the components $R_i$ must be dominant,
        which proves $(4)$. Similarly for $(5)$ working instead with the Kontsevich moduli
        of curves $\cM_{g,1}(X)\to \cM_{g,0}(X)$ the result follows. For $(6)$, the proof
        for rationally connected varieties is contained in \cite[$4.4.(5)$]{debarre}. Let
        $\cC\to U$ be a family which makes $X$ elliptically connected and let $X'\to X$ be
        a connected finite \'etale cover. Consider the pullback diagram and $\cC'\to U'\to
        U$ the Stein factorisation
        \eqn{
            \xymatrix{
                & \cC'=\cC\times_X X'\ar[dl]\ar[r]\ar[d] & X'\ar[d] \\
                U'\ar[dr] & \cC\ar[d]\ar[r] & X \\
                & U. &
            }
        }
        After possibly restricting $U'$ to the open subset of curves in $\cC'$ which are
        irreducible, the family $\cC'\to U'$ makes $X'$ elliptically connected.
    }
    \prop{
        Let $X$ be a smooth projective variety over an algebraically closed field $k$ and
        $f:X\to C$ a flat morphism to a smooth projective curve whose geometric generic
        fibre is separably rationally connected. Then $X$ is $C$-connected.
        \label{rcfibrationovercurveiscconnected}
    }
    \proof{
        From \cite{djs}, there is a section $\s:C\to X$ to $f$. Now from \cite[Theorem 2.13]{kmm92} we can find a
        section to $f$ passing through any two points in different smooth fibres over $C$, hence we can find a copy of
        $C$ passing through two general points. The result now follows from Lemma \ref{gconproperties} parts $(4)$ and
        $(5)$ above after possibly passing to an uncountable extension $K/k$.
    }

    We now come to the main theme of this paper, which is that varieties covered by higher genus curves in a strong
    sense are also covered by rational curves. This is illustrated in the following proposition, and continues in
    the next sections.

    \prop{
        Let $X$ be a $C$-connected variety of dimension at least $3$ over an algebraically
        closed field $k$. Then $X$ is uniruled.
        \label{cconisuniruled}
    }
    \proof{
        We may assume $X$ is projective. Let $u:C\times Y\to X$ be a family such that $u^{(2)}: C\times C\times Y\to
        X\times X$ is dominant. We have $\dim Y + 2 \geq 2\dim X$ and so if $\dim X\geq 3$ we obtain $\dim Y \geq 4$.
        Now, pick general points $x\in X, c\in C$ and denote by $Z\subset Y$ the locus of curves $u_z: C_z\to X$ such
        that $x = u_z(c)$ for all $z\in Z$. We have that $\dim Z \geq \dim Y-(\dim X-1)-1$ and so for $\dim X\geq3$,
        $\dim Z\geq 1$. Since any two general points in $X$ can be connected by the image of a $C_y$, it follows that
        $Z$ does not get contracted to a point when mapped to $\Hom(C,X;c\mapsto x)$. From Bend and Break (see
        \cite[Prop. $3.1$]{debarre}) we obtain a rational curve through $x$ and hence through every general point.
        After possibly an extension to an uncountable algebraically closed field this implies that $X$ is uniruled (see
        \cite[Remark $4.2(5)$]{debarre}).
    }

    If $C$ has genus one, the above result is also proved in Section \ref{sectsix}, even allowing the curve $C$ to vary
    in moduli and with the dimension of $X$ assumed greater or equal to two. On the other hand, a $C$-connected surface
    does not have to be uniruled when $C$ has genus at least two. Consider $C\subset A$ a curve in an abelian surface
    such that $C$ contains the identity $0$ of $A$ and the genus of $C$ is at least two. Consider the map $\phi:C\times
    C\to X$ sending $(p,q)$ to $p-q$. If the image is one dimensional, it has to be isomorphic to $C$ since it has to be
    irreducible and contains the image of $C\times\{0\}$. On the other hand, the image will be closed under the group
    operation, hence would have to be abelian itself, which is a contradiction. Hence $\phi$ is surjective, and we
    obtain that for any $x\in A$, there is a $(p,q)\mapsto x$, hence a morphism $C\cong C\times\{q\}\to X$ passing
    through $x$ and $0$ (for $(q,q)$). Take any two points $x,y\in A$, and consider the image of a morphism from $C$
    through $0$ and the point $x-y$ that we just constructed. Translate this curve by $y$ and obtain an image of $C$
    through $x,y$.

    Denoting by $X\dashrightarrow R(X)$ the maximal rationally chain connected (MRC)
    fibration, we let $R^0(X)=X$, $R^i(X) = R(R^{i-1}X)$ and obtain a tower of MRC
    fibrations 
    \eqn{
        X\dashrightarrow R^1(X) \dashrightarrow \cdots \dashrightarrow R^n(X).
    }
    This tower eventually stabilises, and if $R^i(X)$ is uniruled then $\dim R^{i+1}(X) <
    \dim R^i(X)$. In characteristic zero, we in fact have $R(X) = \ldots = R^n(X)$ (see
    discussion below). In positive characteristic it can be that the tower has length
    greater than one - see the example given in the last section of this paper.

    \prop{
        Let $X$ be a normal and proper $C$-connected variety over an algebraically closed
        field where $C$ is a smooth projective curve. Then the tower $X \dashrightarrow
        R^1(X)\dashrightarrow \cdots \dashrightarrow R^n(X)$ of MRC quotients terminates
        in either a point, a curve or a surface.
        \label{cconmrctopointiorcurve}
    }
    \proof{
        Let $C\times Y\to X$ be a family which makes $X$ a $C$-connected variety. From
        Lemma \ref{gconproperties} part $(1)$ it follows that $R^i(X)$ are also
        $C$-connected. From Proposition \ref{cconisuniruled} we obtain that $R^i(X)$ is
        uniruled if $\dim R^i(X) \geq3$. This implies that $R^{i+1}(X)$ must have
        dimension strictly less than $R^i(X)$ and so the result follows.
    }

    Note that if $k$ is algebraically closed of characteristic zero then we know from \cite{ghs}
    that the MRC quotient $R(X)$ is not uniruled, so if $X$ is $C$-connected of dimension at least three, $R(X)$ must be
    a surface, curve or point, in which case $X$ is respectively a rationally connected fibration over a surface or
    curve, or a point (and so $X$ is rationally connected). From Proposition \ref{rcfibrationovercurveiscconnected} the
    converse holds too for a fibration over a curve. 

    \rem{
        As observed in \cite[Remark $4$]{occhetta}, if the MRC quotient of a smooth complex projective variety $X$ is a
        curve, then the MRC fibration extends to the whole variety and coincides with the Albanese map. 
        \label{occhetta}
    }

\section{Definition of curve connectedness: Free morphisms}

    In this section we define ways in which a morphism from a curve $C$ to a variety $X$
    can deform enough to give a large family of morphisms from $C$ so as to cover $X$. A
    notion studied extensively by Hartshorne \cite{hartamplesubvar} is that of a (local
    complete intersection) subvariety $Y$ in a smooth projective variety $X$ such that the
    normal bundle $\sN_{Y/X}$ is ample. Hartshorne proved in
    \cite[III$.4$]{hartamplesubvar} that for some $g\geq0$ there exists a curve $C\subset
    X$ of genus $g$ such that $\sN_{C/X}$ is ample. Alternatively, Ottem \cite{ottem}
    defines an ample closed subscheme $Y\subset X$ of codimension $r$ to be one where the
    exceptional divisor $\sO(E)$ of the blowup $\Bl_Y X$ of $X$ along $Y$ is an
    $(r-1)$-ample line bundle in the sense that for every coherent sheaf $\sF$ there is an
    integer $m_0>0$ such that $H^i(X, \sF\otimes \sO(E)^m)=0$ for all $m>m_0$ and $i>r-1$.
    One can then prove that if $Y$ is a local complete intersection subscheme of $X$ which
    is ample, then the normal bundle $\sN_{Y/X}$ is an ample bundle. We impose the
    following stronger positivity condition.

    \defn{(\cite[II$.3.1$]{kollar})
        Let $C$ be a smooth proper curve and $X$ a smooth variety over a field $k$. Let
        $f:C\arrow X$ a morphism and $B\subset C$ a closed subscheme with ideal sheaf
        $\sI_B$ and $g=f|_B$. The morphism $f$ is called \textit{free over $g$} if it is
        non-constant and one of the following two equivalent conditions is satisfied:
        \begin{enumerate}
            \item for every $p\in C$ we have $H^1(C, f^*\sT_X\otimes \sI_B(-p))=0$ or,
            \item $H^1(C, f^*\sT_X\otimes \sI_B)=0$ and $f^*\sT_X\otimes \sI_B$ is generated by
                global sections.
        \end{enumerate}
        \label{freeoverBdefn}
    }
    Note that there is also a relative version of the above definition discussed in \cite{ksct}. 
    
    \defn{
        We say that a curve $f:C\arrow X$ is \textit{$r$-free} if for all effective
        divisors $D$ of degree $r\geq0$, $H^1(C, f^*\sT_X\otimes \sO_C(-D))=0$ and
        $f^*\sT_X \otimes \sO_C(-D)$ is generated by global sections. A $0$-free curve is
        called \textit{free} whereas a $1$-free curve is called \textit{very free}. 
        \label{rfreecurvedefn}
    }

    The condition of $r$-freeness makes formal the notion that the curve $C$ deforms in
    $X$ while keeping any general $r$ points fixed. The following follows immediately from
    Lemma \ref{appendixvfimpliesf}.

    \lem{
        If $f:C\to X$ is an $r$-free curve then $f$ is $r'$-free for all $r'\leq r$.
        \label{vfimpliesf}
    }
    
    In the case of $C=\bP^1$, $f^*\sT_X = \oplus_{i=1}^n \sO_{\bP^1}(a_i)$ with $a_1\leq
    \ldots\leq a_n$ so it follows that $f:\bP^1\to X$ is $r$-free if and only if $a_1\geq
    r$. 

    \rem{
        We should remark at this point that there do not exist complete intersection curves of large enough degree which
        are free on a general smooth hypersurface. For example, let $X$ be a degree $d$ smooth hypersurface in $\bP^n$.
        Assume $d\leq n$ since otherwise $X$ will be of general type or Calabi-Yau and will not have any free curves.
        Let $Y_i$ be $n-2$ suitably general hypersurfaces in $\bP^n$ all of degree $e$ and let $C = X \cap_{i=1}^{n-2}
        Y_i$ be the resulting curve. The degree of $C$ is $de^{n-2}$ and the normal bundle is 
        \eqn{
            \sN_{C/X}=\oplus_{i=1}^{n-2} \sO_{\bP^n}(Y_i)|_C = \oplus_{i=1}^{n-2}
            \sO_{\bP^n}(e)|_C. 
        }
        By adjunction, we compute
        \eqn{
            \deg \sT_C = -\deg \omega_C = -d(-n-1+d+\sum_{i=1}^{n-2} e).
        }
        Even setting $e=1$ to make $\deg \sT_C$ as large as possible, and taking into account that $\deg \sN_{C/X} =
        e(n-2)$, we see that $\deg \sT_X|_C = \deg \sT_C + \deg \sN_{C/X}$ is not going to be positive for large values
        of $d$ and $n$.  Positivity of the degree of $\sT_X|_C$ would be necessary for any ampleness conditions. See
        \cite{gounelas} for a discussion on separable rational connectedness of Fano complete intersections.
    }

    A result of Koll\'ar (\cite[II$.1.8$]{kollar}) implies that if the dimension of $X$ is at
    least $3$, a general deformation of a $2$-free morphism is an embedding into $X$. We
    will see (Theorem \ref{rfreeiffample}) that if the genus of $C$ is at least one, this
    holds for any free morphism too. From \cite[II$.3.2$]{kollar}, if a family of curves
    mapping to a variety has a member which is free over $g$, then the locus of all such
    curves in this family is open.

    \lem{
        Let $X$ be a smooth variety over an algebraically closed field $k$, $D\subset X$ a
        divisor and $f:C\to X$ a free morphism. If $p\in C$ then there exists a
        deformation $f':C\to X$ with $f'(p)\notin D$.
        \label{defmissingadivisor}
    }
    \proof{
        By semicontinuity let $U\subset\Hom(C, X)$ be a connected open neighbourhood of
        $[f]$ such that $H^1(C, f_t^*\sT_X) = 0$ for all $[f_t]\in U$. From \cite{mori} it
        follows that the dimension of $U$ is $h^0(C, f^*\sT_X)$. Denote by $\sI_p$ the
        ideal sheaf on $C$ of the closed subscheme with unique point $p$. Since $f$ is
        free, we have $H^1(C, f_t^*\sT_X\otimes \sI_p) = 0$ for all $[f_t]\in U$ and so by
        fixing a point $x\in X$ such that $p\mapsto x$, we have
        \eqn{
            \dim(\Hom(C,X;p\mapsto x) \cap U) &=& h^0(C, f^*\sT_X\otimes \sI_p)\\ &=&
            h^0(C, f^*\sT_X) - \dim X \\ &=& \dim U - \dim X.
        }
        Next, denote by
        \eqn{
            V = \{[f_t] \in U \mid f_t(p)\in D\} = \bigcup_{x\in D}\{[f_t]\in U \mid f_t(p)=x\}
        }
        the subspace of all morphisms in $U$ which send $p$ to a point in the divisor $D$. It
        follows that
        \eqn{
            \codim(V, U) &\geq& \dim U - \dim V \\ &=& h^0(C,f^*\sT_X) - (h^0(C, f^*\sT_X) -
            \dim X + \dim X - 1) = 1
        }
        and hence there exists an $[f']\in U\setminus V$ such that $f'(p)\notin D$.
    }

    \prop{
        Let $X$ be a smooth variety over an algebraically closed field $k$ and $f:C\to X$
        a smooth projective curve which is free over $B\subset C$ a closed subscheme with
        ideal sheaf $\sI_B$. Let $g:X\dashrightarrow Y$ be a generically smooth dominant
        rational map to a smooth proper variety $Y$. Then it follows that $f':=g\circ f:
        C\dashrightarrow Y$ can be deformed to a morphism free over $B$. 
        \label{vfpushesdown}
    }
    \proof{
        Deform $f:C\to X$ so that it misses the codimension $2$ exceptional locus of $g$
        (from \cite[II$.3.7$]{kollar}) so we can assume that the composition $g\circ f:
        C\dashrightarrow Y$ is in fact a non-constant morphism. Starting with the standard exact
        sequence of tangent bundles on $X$ and applying $f^*$ and tensoring with $\sI_B$
        we obtain 
        \eqnn{
            0\to f^*\sT_{X/Y}\otimes \sI_B\to f^*\sT_X\otimes \sI_B \to (g\circ
            f)^*\sT_Y\otimes \sI_B.
            \label{sse1}
        }
        From \cite[Ex. $6.2.10$]{liu} this is exact on the right and we conclude.
    }

    In the case of higher genus curves there exist genus $g$ connected varieties which do
    not have a free or very free curve for all $g\geq1$, for example consider
    $E\times\bP^1$ where $E$ is an elliptic curve. As pointed out after Definition
    \ref{defncconn}, $E\times\bP^1$ is $E$-connected yet it is not possible
    that there exists a morphism $f:C\to E\times\bP^1$ from a curve $C$ such that
    $f^*\sT_{E\times\bP^1}$ is ample since this bundle is isomorphic to $\sO_C\oplus
    \sO_C(2)$ which has a non-ample quotient $\sO_C$. One can however prove the following
    proposition.

    \prop{
        Let $X$ be a smooth variety over an algebraically closed field and $f:C\arrow X$
        a very free morphism for some smooth projective curve $C$. Then $X$ is separably
        $C$-connected.
        \label{vfreesogcon}
    }
    \proof{
        Let $[f]\in Y\subset \Hom(C,X)$ be an open and smooth neighbourhood with cycle map
        $u:C\times Y \arrow X$. We first show that the evaluation map
        \eqn{
            \phi^{(2)}(p,q,f): H^0(C, f^*\sT_X)\arrow f^*\sT_X\otimes k(p) \oplus
            f^*\sT_X\otimes k(q)
        }
        is surjective for $p\neq q$ general points in $C$. Consider the following exact
        sequences of sheaves
        \eqn{
            0 \to f^*\sT_X(-p-q) \to f^*\sT_X \to (f^*\sT_X\otimes k(p)) \oplus
            (f^*\sT_X\otimes k(q)) \to 0 \\
            0 \to f^*\sT_X(-p-q) \to f^*\sT_X(-p) \to f^*\sT_X(-p)\otimes k(q) \to 0
        }
        and note that by taking the long exact sequence in cohomology of the first, to
        show that $\phi^{(2)}(p,q,f)$ is surjective, we need to show that $H^1(C,
        f^*\sT_X(-p-q))=0$. Since $f$ is very free we have from the second sequence that
        $H^0(C, f^*\sT_X(-p))\to f^*\sT_X(-p)\otimes k(q)$ is surjective and also that
        $H^1(C, f^*\sT_X(-p))=0$ from which it follows that $H^1(C, f^*\sT_X(-p-q))=0$.
        Since $\phi^{(2)}(p,q,f)$ is surjective, it follows from \cite[II$.3.5$]{kollar}
        that $u^{(2)}:C\times C\times Y\arrow X\times X$ is smooth at $(p,q,[f])$. We
        conclude that $X$ is separably $C$-connected and thus also separably connected by
        genus $g$ curves.
    }

    \rem{
        It follows that in the setting above that a very free curve (or in fact even a $C$
        such that $X$ is $C$-connected) has the property that it intersects non-trivially
        all but a finite number of divisors. This follows from the fact that we can cover
        an open subset by images of $C$, whose complement will be a proper closed subset
        of $X$ and so contains a finite number of divisors.
    }

\section{Proving uniruledness and rational connectedness}

    In this section we prove that the existence of a free curve of genus $g\geq1$ is equivalent to the existence of an
    $r$-free curve of genus $g$ for all $r\geq1$, and that in characteristic zero this is also equivalent to the
    existence of a very free rational curve. This is in stark contrast to rational curves, where uniruled varieties
    (possessing free rational curves) are not always rationally connected (possessing very free rational curves). We
    begin by noting that there is another type of positive curve one can consider for a smooth projective variety $X$,
    namely $f:C\to X$ such that $f^*\sT_X$ is ample. Note that such a curve automatically has $\sN_{C/X}$ ample. Such
    curves have traditionally been studied in terms of foliations (cf. Theorem \ref{bm}). We will also prove that the
    existence of a curve such that $f^*\sT_X$ is ample is in fact equivalent to the existence of a free curve of the
    same genus.

    \prop{
        Let $X$ be a smooth projective variety over an algebraically closed field $k$ and
        $f:C\to X$ a morphism from a smooth projective curve of genus $g$ such that
        $f^*\sT_X$ ample. Then $X$ is uniruled.
        \label{ampleimpliesuniruled}
    }
    \proof{
        The proof follows the usual Mori argument so we present only a sketch (cf. Theorem \ref{bm}). Note that if $X$
        is a curve, then since a bundle is ample if and only if its pullback under a finite morphism is ample, we obtain
        that $X=\bP^1$. In characteristic zero, after spreading out over a finitely generated extension $\Spec S$ of
        $\Spec \bZ$, one can reduce to any closed prime and consider the equivalent set-up in positive characteristic.
        After pulling back by Frobenius, Lemma \ref{ampleisvgincharp} implies that there is a morphism $f_p^{(n)}: C_p
        \to X_p$ such that $(f_p^{(n)})^*\sT_{X_p}$ is very free (or $r$-free even), where $f_p:C_p\to X_p$ the
        reduction of $f:C\to X$. Bend and Break now produces a rational curve passing through a general point, of
        bounded degree independent of $p$ (see \cite[Prop. $3.5$]{debarre}). These are points in fibres over $\Spec S$ of
        a finite type relative moduli $\Hom_S^d(\bP^1_S, \cX/S, s)$, for $s:\Spec S\to \cX$ a section specifying the
        general point the rational curve goes through. Hence by Chevalley's Theorem the generic fibre over $\Spec S$ is
        also non-empty, and there is a rational curve through a general point of $X$.
    }

    \theo{
        Let $X$ be a smooth projective variety over an algebraically closed field $k$ and
        $f:C\to X$ a morphism from a smooth projective curve of genus $g$ such that
        $f^*\sT_X$ is ample. 
        \begin{enumerate}
            \item If the characteristic $p$ of $k$ is zero, then $X$ is rationally
            connected. 
            \item If $p>0$ then the tower of MRC fibrations terminates with a point.
        \end{enumerate}
        \label{mrctowerispoint}
    }
    \proof{
        From \ref{ampleimpliesuniruled}, we conclude that $X$ is uniruled, regardless of the characteristic. Denote by
        $\pi:X\to R(X)$ the MRC fibration ($R(X)$ is defined up to birational transformation so we may assume $\pi$ is a
        morphism). In characteristic zero, the composition $g:C\to X\to R(X)$ again has $g^*\sT_{R(X)}$ ample, since
        from the proof of \ref{vfpushesdown} the quotient of an ample bundle is ample. So by the Graber-Harris-Starr
        Theorem, since $R(X)$ is uniruled by Proposition \ref{ampleimpliesuniruled}, it must be a point. In positive
        characteristic, it may not be the case that the composition $g:X\to R(X)$ is generically smooth, in which case
        $g^*\sT_{R(X)}$ might not be ample. From Lemma \ref{ampleisvgincharp} however there is a morphism $h:C'\to X$
        such that $h^*\sT_X$ is very free (here $C'$ is a Frobenius pullback of $C$ so has the same genus). From
        \ref{vfreesogcon}, $X$ is separably $C'$-connected, and so by \ref{cconmrctopointiorcurve} we do obtain that the
        tower of MRC quotients $X\to R(X) \to \cdots R^n(X)$ ends in a point, curve or surface. If $\pi: X\to T$ where
        $T:=R^n(X)$ is a smooth projective curve, then by Lemma \ref{defmissingadivisor}, for a point $p\in C'$, we can
        deform $h$ so that the image of $p$ misses the inverse image under $\pi$ of $\pi(h(p))$. Hence $\Hom(C',T)$ is
        at least one dimensional and from de Franchis' Theorem \cite[8.27]{acgh} it follows that $T$ has genus zero or
        one. One excludes the case where $C', T$ both of genus one, by using the fact that there are only countably many
        isogenies between two elliptic curves. Also, $T$ cannot be rational since we have assumed the tower is maximal.
        If now $R^n(X) = S$ is a smooth projective surface, we may assume by pulling back by Frobenius from
        \ref{ampleisvgincharp} and deforming, that there is an at least one dimensional family of morphisms sending a
        fixed point on $C$ to a fixed point on $S$. Hence by Bend and Break \cite[Prop. $3.1$]{debarre} the surface
        would have to be uniruled and we are reduced to the case of a point again.
    }

    Assuming ampleness and regularity of a foliation on a smooth curve in characteristic zero, results of this type have
    been demonstrated in the work of various people, starting with Miyaoka's uniruledness criterion \cite[Theorem
    $8.5$]{miyaoka}. A short summary of recent results follows.

    \theo{(\cite[Theorem $0.1$]{bm}, \cite[Theorem $1$]{ksct})
        Let $X$ be a normal complex projective variety and $C\subset X$ a complete curve
        in the smooth locus of $X$. Assume that $\sF\subset \sT_X$ is a foliation regular
        along $C$ and such that $\sF|_C$ is ample. If $x\in C$ is any point, the leaf
        through $x$ is algebraic and if $x\in C$ is general then the closure of the leaf
        is also rationally connected.
        \label{bm}
    }

    Using \cite[Corollary $0.3$]{bdppnew}, Peternell proved a weaker version of Mumford's
    conjecture on numerical characterisation of rationally connected varieties from which
    one can deduce the following theorem.

    \theo{(\cite[$5.4$, $5.5$]{peternell06})
        Let $X/\bC$ be a projective manifold and $C\subset X$ a possibly singular curve. If
        $\sT_X|_C$ is ample then $X$ is rationally connected. If $\sT_X|_C$ is nef and
        $-K_X.C>0$ then $X$ is uniruled.
        \label{peternell}
    }

    The precise relation between $r$-free morphisms and morphisms $f:C\to X$ such that
    $f^*\sT_X$ is ample is given in the following.

    \theo{
        Let $X$ be a smooth projective variety over an algebraically closed field $k$ and $r\geq0$ any integer. Then
        there exists a morphism $f:C\to X$ from a smooth projective genus $g\geq 1$ curve $C$ such that $f^*\sT_X$ is
        ample if and only if there is an $r$-free morphism $h:C'\to X$ from a genus $g$ smooth projective curve $C'$.
        \label{rfreeiffample}
    }
    \proof{
        Assuming the existence of $h$, we obtain from Lemma \ref{appendixvfimpliesf} that $h$ is also free, and so by
        Proposition \ref{ggandh1impliesample}, $h^*\sT_X$ is ample. If $f^*\sT_X$ is ample, one needs to separate
        between characteristic $p>0$ or equal to zero. In the former case, as in the proof of \ref{ampleimpliesuniruled}
        we get $h:C'\to X$ (here again $C'$ is a Frobenius pullback of $C$ so of genus $g$) which is $r$-free.  When the
        characteristic is zero, $X$ will be rationally connected from \ref{mrctowerispoint}. The idea now is to attach
        many very free rational curves to $C$, apply standard smoothing of combs techniques and prove that the resulting
        general smooth deformations of the comb will be $r$-free genus $g$ curves (cf \cite[II.$7.10$]{kollar}). This
        proceeds as follows. Assemble a comb $D = C\cup \cup_{i=1}^m C_i$ with $m$ rational teeth that are $(r+1)$-free
        like in \cite[II.$7$]{kollar}. For $m$ large enough, $D$ is smoothable to a flat proper family $Y\to T$ where
        the general fibre is isomorphic to $C$, the central fibre is a subcomb of $D$ with a large number of teeth
        depending on $C\subset X$ and $m$, and there is a morphism $F:Y\to X$ which extends $D\to X$. To show that the
        general nearby fibre $f_t:Y_t\to X$ is $r$-free, it suffices to show that $H^1(Y_t, f_t^*\sT_X(-\sum_{i=0}^r
        p_i))$ for $p_0, p_1,\ldots, p_r$ any points on $Y_t\subset X$ (see Definition \ref{freeoverBdefn}). Pick
        sections $s_0,s_1,\ldots,s_r: T\to Y$ with $s_i(t)=p_i$. Let $E=F^*\sT_X(-\sum_{i=1}^rs_i(T))$. By Riemann-Roch,
        for $m$ large enough, we have that $H^1(C, M\otimes E|_C)=0$ for all line bundles $M$ of degree larger than $m$,
        and also that $E|_{C_i}$ is ample since $C_i$ is $(r+1)$-free. Now apply \cite[II.$7.10.1$]{kollar} for $m$
        large enough.
    }

    Using any of Theorems \ref{bm}, \ref{peternell} or \ref{mrctowerispoint}, a smooth
    projective variety $X$ over an algebraically closed field of characteristic zero with
    a free genus $g$ curve $f:C\to X$ such that $g\geq1$ is automatically rationally
    connected.

    \rem{
        At this point we cannot prove that in positive characteristic, assuming that we have a free curve $f:C\to X$ of
        genus $g\geq 1$ implies that $X$ is separably rationally connected or even rationally chain connected. It is
        tempting to hope that both statements are true though. Jason Starr informs us that his maximal free rational
        quotient (MFRC) \cite{starr06} gives a generically (on the source) smooth morphism $X\to R_f(X)$ over any
        algebraically closed field $k$, so if $X$ contained a free rational curve $f:\bP^1\to X$, then $\dim R_f(X) <
        \dim X$. Hence, if $f:C\to X$ a free curve of genus $g\geq1$ implied that we have a free rational curve
        $\bP^1\to X$ (we do not know how to show this), taking successive MFRC quotients and using Proposition
        \ref{vfpushesdown} would reduce the tower of MFRC quotients to a point. This does not mean that $X$ will
        necessarily be rationally connected, but since there is a free rational curve on $X$, it will at least be
        separably uniruled. Even though Bend and Break arguments give us the existence of many rational curves, the
        author does not know any general techniques to construct free rational curves in positive characteristic. See
        the last two sections for results in this direction.
        \label{conjectureremark}
    }

\section{Elliptically connected varieties}\label{sectsix}

    In this section we will study more carefully the case of genus one. Denoting RC and EC
    to mean rationally and elliptically connected (genus one connected) respectively, we
    have the following inclusions of sets of varieties
    \eqn{
        \{\text{rational}\} \subsetneq \{\text{unirational}\} \subseteq \{\text{RC}\}
        \subsetneq \{\text{EC}\} \subsetneq \{\text{uniruled}\}.
    }
    It is an open problem whether there exists a non-unirational rationally connected variety but it is widely expected
    these do exist. The following result is in the spirit of \ref{cconisuniruled}. The following proof was suggested by
    the anonymous referee.

    \prop{
        Let $X$ be an elliptically connected smooth projective variety of $\dim X\geq 2$
        over an algebraically closed field $k$. Then $X$ is uniruled.
        \label{uniruled}
    }
    \proof{
        Like in \ref{cconisuniruled}, for $\cC\to U$ a family of genus one curves mapping
        to $X$ such that $\cC\times_U\cC \to X\times X$ is dominant, there is an at least
        one dimensional locus $Z\subset U$ parametrising curves which pass through a
        (general) point $x\in X$. In fact, after fixing a general hyperplane $H$, we obtain a morphism
        $Z\to \cM_{1,2}(X)$ where for $z\in Z$, the two marked points are the point
        $p_z\in \cC_z$ sent to $x$, and a point $q_z\in\cC_z$ which is sent to $H$.
        Denote also by $\cC\to Z$ the restriction of the family from $U$. Consider now a
        compactification and the induced rational map to $X$
        \eqn{
            \xymatrix{
                \bar{\cC}\ar[d]^\pi\ar@{-->}[r]^f & X \\
                \bar{Z} & 
            }
        }
        
        and let $\mu:\bar{Z}\to\bar{\cM}_{1,2}(X)$ be the moduli map. Since $\cM_{1,2}$
        contains no proper subvarieties which do not get contracted when mapped to
        $\cM_1$, either the image of $\mu$ meets the boundary, which implies that there is
        a rational curve through $x$, or $\mu$ is a contraction to a point. In the latter
        case, we thus have that the family $\pi$ is isotrivial, so after passing to a
        finite flat cover $\bar{Z}'$ of $\bar{Z}$ we obtain $C\times\bar{Z}'\to\bar{Z}'$,
        with $f':C\times\bar{Z}' \dashrightarrow X$ the induced morphism. From the
        construction, we also obtain a point $p\in C$ (mapped to each $p_z$ under the map
        $C\times\bar{Z}'\to\bar{\cC}$) such that $f'$ contracts $\{p\}\times\bar{Z}'$ to
        $x$. If $f'$ were defined everywhere, Mumford's Rigidity Theorem would imply that
        all fibres $\{s\}\times\bar{Z}'$ are contracted, which contradicts the fact that
        images of our initial family dominate $H$. Hence $f'$ is not defined everywhere
        and like in Bend and Break, we obtain a rational curve through $x$.
    }

    \theo{
        Let $X$ be a smooth projective variety over an algebraically closed field $k$ of
        characteristic zero. Then $X$ is elliptically connected if and only if it is
        rationally connected or a rationally connected fibration over an elliptic curve.
        \label{elconclassification}
    }
    \proof{
        Consider the MRC fibration $\pi: X \dashrightarrow R(X)$ where $R(X)$ is
        elliptically connected as $\pi$ is dominant. Since $R(X)$ is elliptically
        connected and not uniruled, it follows from Proposition \ref{uniruled} that it
        must be either of dimension $0$ and thus $X$ is rationally connected, or of
        dimension $1$ and so an elliptic curve $E$ by Riemann-Hurwitz. By Remark
        \ref{occhetta}, the MRC fibration coincides with the map to the Albanese and so
        fibres of $X\to E'$ are rationally connected. Conversely, we have seen that a
        rationally connected variety is elliptically connected in Lemma
        \ref{genusgconnisgconforbiggerg}. If on the other hand $X$ is a rationally
        connected fibration over an elliptic curve $E$ then from Proposition
        \ref{rcfibrationovercurveiscconnected} we know that it is $E$-connected.
    }

    If $k$ is of positive characteristic, using the same methods as in Lemma
    \ref{cconmrctopointiorcurve} we deduce that for an elliptically connected variety,
    the tower of MRC fibrations terminates with a point or a curve. 

    \rem{
        Note that Bjorn Poonen \cite{poonen} has constructed non-trivial examples over an
        arbitrary field, of elliptically connected threefolds which are not rationally
        connected. These are Ch\^atelet surface fibrations over an elliptic curve.
    }

\section{Towards a positive characteristic analogue}

    From Remark \ref{conjectureremark} and the work preceding it, we would like to
    demonstrate that the existence of a free higher genus curve implies the existence of a
    free rational curve in positive characteristic, something which holds in
    characteristic zero from Theorem \ref{mrctowerispoint}. In this section we make
    the first steps in this direction. If $f:C\to X$ is a very free morphism from a smooth
    projective curve of genus $g\geq2$ to a smooth projective variety $X$, then $K_X.C =
    -\deg f^*T_X < 0$ from the ampleness of $f^*T_X$. In fact, a Riemann-Roch calculation
    gives a better bound of $K_X.C\leq -n(g-1)$ where $n = \dim X$.

    \prop{
        Let $X$ be a smooth projective surface over an algebraically closed field $k$ with
        $f:C\to X$ a free morphism from a smooth projective curve $C$ of genus $g>0$ or a
        very free morphism of genus zero. It follows that $X$ is separably rationally
        connected.
    }
    \proof{
        If $C$ is of genus zero then $X$ is separably rationally connected by definition.  From the discussion above we
        have that $K_X$ is not nef. Also, any surface $Y$ which is birational to $X$ admits a morphism $C\to Y$ from
        \ref{vfpushesdown}, which is again free, so $K_Y$ is also not nef. From the classification of surfaces this
        means that $X$ is either rational or ruled. If ruled, $X$ would admit a birational morphism to $\bP^1\times C$.
        The free morphism $f:C\to X$ would give a free morphism $C\to \bP^1\times C$ which would mean $C$ is $\bP^1$ and
        $X$ was rational.
    }
    \rem{
        Some remarks about the case of dimension three, where the minimal model program is
        incomplete in positive characteristic. From the main theorem in \cite{kollar91},
        assuming $X$ is smooth and that it admits a free morphism from a curve, we can
        contract extremal rays in the cone of curves in arbitrary characteristic, to
        obtain a Fano fibration over a curve, surface or point. In the case where there
        exists a conic fibration $X\to S$ where $S$ is a smooth surface, Koll\'ar proves
        that if the characteristic of $k$ is not $2$ then the general fibre is smooth.
        From Proposition \ref{vfpushesdown} it follows that the composition morphism $C\to
        S$ is free and so from the above proposition for the case of surfaces, $S$ is a
        rational surface. Hence $X$ is a conic bundle over a rational surface hence
        separably rationally connected. If $X\to Y$ a Fano fibration over a curve, to the
        author's knowledge, it is not known whether the fibres of the del Pezzo surface
        fibration over $Y$ obtained in this way must be smooth. Assuming for the time
        being that they were, they would be separably rationally connected and from the
        deformation theory argument in Theorem \ref{mrctowerispoint} and de Franchis'
        Theorem \cite[8.27]{acgh}, $Y$ would be $\bP^1$. From the de
        Jong-Starr Theorem we would obtain sections $\bP^1\to X$ from which we could
        assemble combs with very free teeth to be smoothed to very free rational curves in
        $X$, showing that $X$ is separably rationally connected. Finally, even though it
        is open whether Fano threefolds are separably rationally connected (this result is
        not true in higher dimensions however), Shepherd-Barron \cite{sb} proved that Fano
        threefolds of Picard rank one are liftable to characteristic zero, hence admitting
        a very free morphism implies they are separably rationally connected.
    }

    The following result is well known in the case of $\bP^1$ (see \cite[$4.18$]{debarre})
    and easily extends to higher genus.

    \prop{
        \label{freecurveimpliesnoglobalforms}
        Let $f:C\to X$ be a very free morphism from a smooth projective curve $C$ to a
        smooth projective variety $X$ over an algebraically closed field $k$. Then for all
        positive integers $m,\ell$
        \eqn{
            H^0(X, (\Omega^\ell_X)^{\otimes m})=0.
        }
    }
    \proof{
        Since $f:C\to X$ is very free, from Proposition \ref{vfreesogcon} there is a
        variety $U$ such that $C\times U\to X$ makes $X$ separably $C$-connected. Being
        very free is an open property (\cite[II$.3.2$]{kollar}) so we can assume that the
        general morphism $f_u:C_u\to X$ for $u\in U$ is very free and also an immersion
        from \cite[II.$1.8$]{kollar}, and so $f_u^*\sT_X$ is ample from Proposition
        \ref{ggandh1impliesample} (and by definition of a very free curve in the genus zero 
        case). We conclude that for a general point $x\in X$ there is a morphism
        $f_u:C_u\to X$ such that $f_u^*\sT_X$ is ample and whose image passes through $x$.
        Hence since $f_u^*\Omega^1_X$ is negative, any section of $(\Omega_X^\ell)^{\otimes
        m}$ must vanish on the image $f(C_u)$ hence on a dense open subset of $X$, and so
        on $X$.
    }

    \cor{
        Let $f:C\to X$ as above. Then the Albanese variety $\Alb X$ is trivial.
    }
    \proof{
        Note that we have that $\dim\Alb X \leq \dim H^1(X, \sO_X) = h^{0,1}$. In
        characteristic zero Hodge duality gives that $h^{1,0}=h^{0,1}$ but more generally
        over any algebraically closed field we have (see \cite{igusa}) that $\dim\Alb X
        \leq h^{1,0} = h^0(X, \Omega^1_X)$. The result follows from Proposition
        \ref{freecurveimpliesnoglobalforms}.
    }

    The above also follows from the result in \cite{gounelas}, which says that in the above situation $H^1(X, \sO_X)=0$.
    See ibid.\ for a discussion around the vanishing of $H^i(X,\sO_X)$ for separably rationally connected varieties
    in positive characteristic. Note also that if $X$ is $C$-connected, since any map $C\to \Alb X$ must factor through
    the Jacobian, and there are only countably many homomorphisms between abelian varieties, one concludes that the
    image of $X$ in $\Alb X$ is either a point or a curve.

\section{An example in positive characteristic}

    Let $X$ be the Fermat quintic surface $x_0^5+x_1^5+x_2^5+x_3^5=0$ in $\bP^3$ over an
    algebraically closed field of characteristic $p$. In \cite{shioda} it is proven that
    if $p\neq 5$ and $p$ is not congruent to $1$ modulo $5$, then $X$ is a unirational
    general type surface and if we quotient by the action of the group $G$ of $5$-th
    roots of unity $x_i\mapsto \zeta^i x_i$, then we obtain a Godeaux surface which is
    again unirational but has algebraic fundamental group $\pi_1^{\et}(X/G, \bar{y})
    \cong \bZ/5\bZ$. Note that in characteristic zero, the notions of rationally chain
    connected, rationally connected, freely rationally connected (see \cite{shen}) and
    separably rationally connected all coincide and it is known that each variety in
    this class is simply connected. In positive characteristic however these notions are
    in decreasing generality and can differ. A rationally chain connected variety always
    has finite fundamental group (see \cite{cl03}) whereas a freely rationally connected
    variety is simply connected (see \cite{shen}). Note that Shioda's example above gives
    a unirational and hence rationally connected variety over a characteristic $p$
    algebraically closed field which is not simply connected.  

    \n
    We show there is a smooth projective variety in characteristic $p$ which has infinite
    \'etale fundamental group but after a finite number of MRC quotients we terminate with
    a point. Let $C$ be a smooth $5$ to $1$ cover of $\bP^1$, with defining affine
    equation of the form $y^5=f(x)$ where $f$ is a general polynomial of high degree. We
    have an action of $G=\bZ/5\bZ$ on $C$ which we can extend to the product $X\times C$
    of the above Fermat quintic $X$ with $C$. Projecting from the quotient onto the second
    factor we have a morphism $(X\times C)/G \to \bP^1$ where we have identified $C/G$
    with $\bP^1$. The general fibre of this morphism is isomorphic to $X$. We have a short
    exact sequence
    \begin{equation*}
        1\to \pi_1^{\et}(X,\bar{x})\times\pi_1^{\et}(C,\bar{c}) \to
        \pi_1^{\et}((X\times C)/G,\bar{z}) \to G\to1.
    \end{equation*}
    Hence we have constructed an example of a smooth projective variety over an
    algebraically closed field of characteristic $p$ whose fundamental group is infinite
    yet whose tower of MRC quotients terminates with a point.

\end{document}